\documentclass{amsart}
\usepackage{appendix}
\usepackage{tabularx}
\usepackage[numbers,sort&compress]{natbib}
\usepackage{color}
\usepackage[colorlinks=true,linkcolor=blue,citecolor=blue,urlcolor=blue]{hyperref}
\newtheorem{theorem}{Theorem}[section]
\newtheorem{lemma}[theorem]{Lemma}
\newtheorem{corollary}[theorem]{Corollary}

\theoremstyle{definition}

\theoremstyle{remark}
\newtheorem{remark}[theorem]{Remark}
\numberwithin{equation}{section}
\allowdisplaybreaks[4]

\begin{document}

\title[$C^{0}$ stability of the Laplacian flow]{$C^{0}$ stability of torsion-free $G_{2}$ structures along the Laplacian flow}


\author{Chuanhuan Li}
\address{ Shanghai Institute for Mathematics and Interdisciplinary Sciences (SIMIS), Shanghai 200433, China \newline
${\quad}$ Research Institute of Intelligent Complex Systems, Fudan University, Shanghai 200433, China}
\curraddr{}
\email{chli@simis.cn}
\thanks{}

\author{Yi Li}
\address{Center for Mathematics and Interdisciplinary Sciences, Fudan University, Shanghai 200433, China \newline
${\quad}$ Shanghai Institute for Mathematics and Interdisciplinary Sciences (SIMIS), Shanghai 200433, China}
\curraddr{}
\email{yilicms@simis.cn, yilicms@gmail.com}



\date{}

\dedicatory{}

\begin{abstract}
In this paper, we prove that torsion-free $G_{2}$-structures are (weakly) dynamically stable along the Laplacian flow for closed $G_{2}$-structures in $C^{0}$-open neighbourhood. This weakens the assumption in \cite{LW19} and confirms Bryant's expectation \cite{Bry06} in the $C^{0}$ norm.

\end{abstract}

\maketitle

\section{Introduction}

In order to study a smooth $7$-manifold $M$ admitting a {torsion-free} $G_{2}$-structure, Bryant \cite{Bry06} introduced the following Laplacian $G_2$ flow for closed $ G_{2}$-structure:
\begin{equation}
  \left \{
       \begin{array}{rl}
          \partial_{ t}\varphi(t)&=\Delta_{\varphi(t)}\varphi(t),\\
           d\varphi(t)&=0,\\
           \varphi(0)&=\varphi_{0},
       \end{array}
  \right.
  \label{eq:flow}
\end{equation}
where $\Delta_{\varphi(t)}=dd^{\ast}_{\varphi(t)}+d^{\ast}_{\varphi(t)}d$ is the Hodge Laplacian of $g(t)$ and $\varphi_{0}$ is an initial  closed $G_{2}$-structure({a positive closed three-form).} Here $g(t)$ is the Riemannian metric {algebraically determined by} $\varphi(t)$. Since $\Delta_{\varphi}\varphi=dd^{\ast}_{\varphi}\varphi$ for a closed $G_{2}$-structure $\varphi$, we see that the closedness of $\varphi(t)$ is preserved along the Laplacian flow $\eqref{eq:flow}$. 

The flow $\eqref{eq:flow}$ can be viewed as the gradient flow for the Hitchin functional when the variations are restricted to the cohomology class of the closed $G_{2}$-structure, where Hitchin functional was introduced by Hitchin \cite{Hit00}
$$\mathcal{H}:[\overline{\varphi}]_{+}\longrightarrow\mathbb{R}^{+},\ \varphi\longmapsto\frac{1}{7}\int_{M}\varphi\wedge\psi=\int_{M}\ast_{\varphi}1.$$
Here $\overline{\varphi}$ is a closed $G_{2}$-structure on $M$ and  $[\overline{\varphi}]_{+}$ is the open subset of the cohomology class $[\overline{\varphi}]$ consisting of $G_{2}$-structures. Any critical point of $\mathcal{H}$ gives rise to a torsion-free $G_{2}$-structure {when $7$-manifold $M$ is compact.} 

In \cite{BX11}, applying De Turck's trick and Hamilton's Nash-Moser inverse function theorem, Bryant and Xu modified the Laplacian flow \eqref{eq:flow} via a term $\mathcal{L}_{V(\varphi(t))}\varphi(t)$ for some vector field $V(\varphi(t))$ and considered the following Laplacian-DeTurck flow
\begin{equation}
  \left \{
       \begin{array}{rl}
          \partial_{ t}\varphi(t)&=\Delta_{\varphi(t)}\varphi(t)+\mathcal{L}_{V(\varphi(t))}\varphi(t),\\
           d\varphi(t)&=0,\\
           \varphi(0)&=\varphi_{0},
       \end{array}
  \right.
  \label{Laplacian-DeTurck flow}
\end{equation}
They proved the short-time existence for $\eqref{eq:flow}$. In 2017, Lotay and Wei \cite{LW17} proved the long-time existence for the Laplacian $G_{2}$ flow $\eqref{eq:flow}$ under bounded Riemann curvature or Ricci curvature. In 2021, the second author \cite{local curvature of G2} gave an alternative proof via local curvature estimates, showing that the flow persists as long as the Ricci curvature remains bounded. The long-time existence of solutions to the Laplacian flow with bounded scalar curvature remains an open problem. Some special cases and partial results can be found in \cite{FY18,LL21,LL26,LLZ26,PS24}.

In \cite{Bry06}, Bryant suggests that it is reasonable to expect the following
\begin{quote}
    {\bf Expectation:} {\it Given a torsion-free $G_{2}$-structure $\bar{\varphi}$, for $\varphi_{0}\in[\bar{\varphi}]_{+}$ sufficiently near $\bar{\varphi}$ in an appropriate norm, the Laplacian flow $\eqref{eq:flow}$ with $\varphi_{0}$ as initial value would converge to a point on ${\rm Diff}^{0}(M)\cdot\bar{\varphi}$.}
\end{quote}
Lotay and Wei \cite{LW19} confirmed this expectation when $\varphi_{0}$  is sufficiently close to $\bar{\varphi}$ in the $C^{9}$-norm. It is therefore natural to consider the optimal norm in which $\varphi_{0}$ should be sufficiently close to $\bar{\varphi}$.

In this paper, we study the optimal norm in which $\varphi_{0}$ should be sufficiently close to $\bar{\varphi}$ in order to confirm Bryant's expectation. We show that closeness in the $C^{0}$-norm is enough. Motivated by the heat semigroup theory in partial differential equations, we derive an estimate for the norm of the perturbed form to $\bar{\varphi}$.
\begin{theorem}\label{thm1.1}
     Let $\bar{\varphi}$ be a torsion-free $G_{2}$ structure on a compact $7$-dimensional manifold $M$. And $\tilde{\varphi}(t)$ be a solution to the Laplacian-DeTurck flow \eqref{Laplacian-DeTurck flow} with initial value $\varphi_{0}$ on $[0,T)$. If $\|\varphi_{0}-\bar{\varphi}\|_{C^{0}_{\tilde{g}}}\leq \delta$, then $T> 1$ and
     $$\sum_{i=0}^{n}\tau^{\frac{i}{2}}\|\bar{\nabla}^{i}(\tilde{\varphi}(\tau)-\bar{\varphi})\|_{C^{0}_{\tilde{g}}}\leq C\delta,\quad \text{for}\quad 0<\tau\leq 1,$$
      where $C$ is constant depending only on $M,\bar{\varphi},n$ and $\alpha$ with $0<\alpha<1$.
\end{theorem}
Combining this with Theorem 5.1 and Page 521-524 in \cite{LW19}, we confirm Bryant's expectation when $\varphi_{0}$  is sufficiently close to $\bar{\varphi}$ in the $C^{0}$-norm.
\begin{theorem}\label{thm1.2}
    Let $\bar{\varphi}$ be a torsion-free $G_{2}$ structure on a compact $7$-dimensional manifold $M$. There is a $C^{0}$-neighbourhood $\mathcal{U}$ of $\bar{\varphi}$ such that for any $\varphi_{0}\in[\bar{\varphi}]_{+}\cap\mathcal{U}$, the Laplacian flow \eqref{eq:flow} with initial value $\varphi_{0}$ exists for all $t\in[0,\infty)$ and converges smoothly to $\varphi_{\infty}\in {\rm Diff}^{0}(M)\cdot\bar{\varphi}$ as $t\rightarrow\infty.$
\end{theorem}
If we assume Joyce's existence result for torsion-free $G_{2}$ structures in \cite{Joy96}, we obtain the following immediate corollary of Theorem \ref{thm1.2}.
\begin{corollary}
    Let $\varphi_{0}$ be a closed $G_{2}$-structure on a compact $7$-manifold $M$. There exists a $C^{0}$-open neighbourhood $\mathcal{U}$ of $0$ in $\Omega^{3}(M)$ such that if $d^{\ast}_{\varphi_{0}}\varphi_{0}=d^{\ast}_{\varphi_{0}}\gamma$ for some $\gamma\in\mathcal{U}$, then the Laplacian flow \eqref{eq:flow} with initial value $\varphi_{0}$ exists for all time and converges to a torsion-free $G_{2}$-structure.
\end{corollary}

We give an outline of this paper. We review the basic theory in Section \ref{section2} about $G_{2}$-structure, $G_{2}$-decompositions of $2$-forms and $3$-forms and define some norms of forms. Section \ref{section3} derives the evolution equation of the perturbed form $\theta(t)$ and obtains some estimates of the heat equation of forms. In Section \ref{section4}, using the estimates in Section \ref{section3}, we get the main estimate of $\theta(t)$ and prove the main result.

${}$

{\bf Declaration of AI}. AI tools were used to assist with language editing, mathematical exposition, and the checking of parts of the arguments. The authors take full responsibility for the final manuscript, including the verification of all mathematical arguments and references.

\section{$G_{2}$-structures}\label{section2}
In this section, we give some basis theory of $G_{2}$-structures in \cite{Bry06,Spiros introduce to G2,LW17}.


Let $\{e_{1},e_{2},\cdots,e_{7}\}$ denote the standard basis of $\mathbb{R}^{7}$ and let $\{e^{1},e^{2},\cdots,e^{7}\}$ be its dual basis. Define the 3-form
$$\phi:=e^{123}+e^{145}+e^{167}+e^{246}-e^{257}-e^{347}-e^{356},$$
where $e^{ijk}:=e^{i}\wedge e^{j}\wedge e^{k}$. The subgroup of ${\rm GL}(7,\mathbb{R})$ {that fixes $\phi$} is the exceptional Lie group $G_{2}$, a compact, connected, simple $14$-dimensional Lie subgroup of ${\rm  SO}(7)$. In fact, $G_{2}$ acts irreducibly on $\mathbb{R}^{7}$ and preserves the metric and orientation for which $\{e_{1},e_{2},\cdots,e_{7}\}$ is an oriented orthonormal basis. Note that $G_{2}$ also preserves the $4$-form
$$\ast_{\phi}\phi=e^{4567}+e^{2367}+e^{2345}+e^{1357}-e^{1346}-e^{1256}-e^{1247},$$
where $\ast_{\phi}$ is the Hodge star operator determined by the metric and orientation.

For a smooth $7$-manifold $M$ and a point $x\in M$, we define as in \cite{local curvature of G2,LW17}
\begin{align}
    \wedge_{+}^{3}(T_{x}^{\ast}M):=\left\{\varphi_{x}\in\wedge^{3}(T_{x}^{\ast}M)\Big{|} u^{\ast}\phi=\varphi_{x},
    \text{for\ invertible}\ u\in\text{Hom}_{\mathbb{R}}(T_{x}M,\mathbb{R}^{7})\right\}\notag
\end{align}
and the bundle
\begin{align}
    \wedge_{+}^{3}(T^{\ast}M):=\bigcup_{x\in M}\wedge_{+}^{3}(T_{x}^{\ast}M)\notag.
\end{align}
We call a section $\varphi$ of $\wedge_{+}^{3}(T^{\ast}M)$ a {\it positive $3$-form} on $M$ or a {\it $G_{2}$-structure} on $M$, and denote the space of positive $3$-forms by $\Omega^{3}_{+}(M)$. The existence of $G_{2}$-structure is equivalent to the
property that $M$ is both orientable and spinnable.

For a {positive} $3$-form $\varphi$, we define a $\Omega^{7}(M)$-valued bilinear form $\text{B}_{\varphi}$ by
$$\text{B}_{\varphi}(u,v)=\frac{1}{6}(u\lrcorner\varphi)\wedge(v\lrcorner\varphi)\wedge\varphi,$$
where $u, v$ are tangent vectors on $M$ and $``\lrcorner"$ is the interior multiplication operator (\textit{Here we use the orientation in \cite{Bry06}}). Then we can see that any $\varphi\in\Omega^{3}_{+}(M)$ determines a Riemannian metric $g_{\varphi}$ and an orientation $d V_{\varphi}$, hence the Hodge star operator $\ast_{\varphi}$ and the associated $4$-form
$$\psi:=\ast_{\varphi}\varphi$$
can also be uniquely determined by $\varphi$.


The group $G_{2}$ acts irreducibly on $\mathbb{R}^{7}$ (and hence on $\wedge^{1}(\mathbb{R}^{7})^{\ast}$ and $\wedge^{6}(\mathbb{R}^{7})^{\ast}$), but it
acts reducibly on $\wedge^{k}(\mathbb{R}^{7})^{\ast}$ for $2\leq k\leq 5$. Hence a $G_{2}$ structure $\varphi$ induces splittings
of the bundles $\wedge^{k}(T^{\ast}M)(2\leq k\leq5)$ into direct summands, which we denote by
$\wedge^{k}_{l}(T^{\ast}M,\varphi)$ with $l$ being the rank of the bundle. We let the space of {smooth} sections
of $\wedge^{k}_{l}(T^{\ast}M,\varphi)$ be $\Omega^{k}_{l}(M)$. Define the natural projections
$$\pi^{k}_{l}:\Omega^{k}(M)\longrightarrow \Omega^{k}_{l}(M),\ \ \alpha\longmapsto \pi^{k}_{l}(\alpha).$$
Then we have
\begin{align}
    \Omega^{2}(M)&=\Omega^{2}_{7}(M)\oplus\Omega^{2}_{14}(M),\notag\\
    \Omega^{3}(M)&=\Omega^{3}_{1}(M)\oplus\Omega^{3}_{7}(M)\oplus\Omega^{3}_{27}(M)\notag.
\end{align}
where each component is determined by
\begin{align}
    \Omega^{2}_{7}(M)&=\{X\lrcorner\varphi:X\in {\Gamma}(TM)\}=\{\beta\in\Omega^{2}(M):\ast_{\varphi}(\varphi\wedge\beta)=2\beta\},\notag\\
    \Omega^{2}_{14}(M)&=\{\beta\in\Omega^{2}(M):\psi\wedge\beta=0\}=\{\beta\in\Omega^{2}(M):\ast_{\varphi}(\varphi\wedge\beta)=-\beta\},\notag
\end{align}
and
\begin{align}
    \Omega^{3}_{1}(M)&=\{f\varphi:f\in C^{\infty}(M)\},\notag\\
    \Omega^{3}_{7}(M)&=\{\ast_{\varphi}(\varphi\wedge\alpha):\alpha\in\Omega^{1}(M)\}=\{X\lrcorner\psi:X\in {\Gamma}(TM)\},\notag\\
    \Omega^{3}_{27}(M)&=\{\eta\in\Omega^{3}(M):\eta\wedge\varphi=\eta\wedge\psi=0\}.\notag
\end{align}

\begin{remark}\label{remark2.2}
    $\Omega^{4}$ and $\Omega^{5}$ have the corresponding decompositions by Hodge duality. The more details for $G_{2}$-decomposition see \cite{Bry06,Spiros introduce to G2}.
\end{remark}

We define the Hodge Laplacian induced by $g_{\varphi}$ as
$$\Delta_{\varphi}=dd^{\ast}_{\varphi}+d^{\ast}_{\varphi}d,$$
and the trace Laplacian induced by $g_{\varphi}$ as
$$\triangle={\rm tr}_{g}(\nabla^{2}).$$
Let $(M,g)$ be a compact Riemannian manifold. For a $k$-tensor $A$, the $C^{0}$-norm is defined by
$$\|A\|_{C^{0}(M,g)}:=\max_{\bar{M}}|A|_{g}.$$
Fix $0<\rho<\operatorname{inj}(M,g)$. The H\"older seminorm compares tensors by $g$-parallel transport $P_{\gamma}$ along the unique short minimizing geodesic $\gamma$ connecting $x$ to $y$:
 $$[A]_\alpha=\sup_{0<d_{g}(x,y)<\rho}
 \frac{|A(x)-P_{\gamma}A(y)|_{g}}{d_{g}(x,y)^\alpha}$$
with $0<\alpha<1$, where $d_{g}(x,y)$ is the distance induced by metric $g$.
We then define the $C^{k}$ norm by
$$
 \|A\|_{C^k(M,g)}=\sum_{i=0}^k\|\nabla^iA\|_{C^0(M,g)},$$
 and the $C^{k,\alpha}$-norm by
 $$
 \|A\|_{C^{k,\alpha}(M,g)}=\|A\|_{C^k(M,g)}+[\nabla^kA]_\alpha.
$$
For convenience, we often suppress the dependence on $(M,g)$ and write 
$\|A\|_{C^k}$  instead of $\|A\|_{C^k(M,g)}$.

\section{estimate for a heat equation}\label{section3}
In this section, we derive the evolution equation for the perturbation of  $\bar{\varphi}$ under Laplacian-DeTurck flow $\eqref{Laplacian-DeTurck flow}$ and obtain some estimates for the heat equation of forms. These will be used in Section \ref{section4}.

Let $M$ be a compact $7$-dimensional manifold admitting a torsion-free $G_{2}$-structure $\bar{\varphi}$. Since the solution $\tilde{\varphi}(t)$ of the  Laplacian-DeTurck flow $\eqref{Laplacian-DeTurck flow}$ lies in the cohomology class $[\bar{\varphi}]$, we can write $$\tilde{\varphi}=\bar{\varphi}+\theta(t),$$ 
where $\theta(t)$ takes values in exact forms.

Now, we recall some definitions of the Laplacian-DeTurck flow in \cite{BX11}. Let $S=\nabla-\bar{\nabla}$ be the tensor taking values in $TM\otimes S^{2}T^{\ast}M$, where $\nabla$ and $\bar{\nabla}$ are  connections with respect to $\varphi(t)$ of the solution to the Laplacian flow $\eqref{eq:flow}$ and the background $\bar{\varphi}$. We let $V(\varphi)$ be
$$V(\varphi)=-5V_{1}(\varphi)-V_{2}(\varphi)$$
as in \cite{BX11}. $V_{1}$ and $V_{2}$ are given locally by
$$V_{1}(\varphi)=\frac{1}{7}g^{pq}S_{pq}^{i}e_{i},\quad V_{2}(\varphi)=\frac{2}{A}(g^{kj}S_{ik}^{i}e_{j}+5V_{1}(\varphi)),$$
where the constant $A$ is from Lemma 1.7 in \cite{BX11}.

Since Lemma 4.1 in \cite{LW19} requires a $C^{k}$ assumption on $\theta$. We therefore also derive the evolution equation for $\theta(t)$ under Laplacian-DeTurck flow $\eqref{Laplacian-DeTurck flow}$.
\begin{lemma}\label{lem3.1}
    For all time $t$ for which the solution $\tilde{\varphi}(t)$ of \eqref{Laplacian-DeTurck flow} exists, we have
    \begin{equation}\label{theta}
        \frac{\partial}{\partial t}\theta(t)=-\Delta_{\bar{\varphi}}\theta(t)+dF(\bar{\varphi},\tilde{\varphi}(t),\theta(t),\bar{\nabla}\theta(t)),
    \end{equation}
    where $F$ is a $2$-form satisfying
    $$F(\bar{\varphi},\tilde{\varphi}(t),\theta(t),\bar{\nabla}\theta(t))=f(\bar{\varphi},\tilde{\varphi}(t))\theta(t)\bar{\nabla}\theta(t)$$
    with some tensor $f$ depending on $\bar{\varphi}$ and $\tilde{\varphi}(t)$. 
    
    In particular, there exists $\delta^{\ast}=\delta^{\ast}(\bar{\varphi})>0$ such that if $\delta\in(0,\delta^{\ast})$ and $\|\theta(t)\|_{C^{0}_{\bar{g}}}<\delta$, then there exists a uniform constant $C$ depending on $\bar{\varphi}$ and $\delta$, such that
    $$|F|_{\bar{g}}\leq C|\theta(t)|_{\bar{g}}|\bar{\nabla}\theta(t)|_{\bar{g}}.$$
\end{lemma}
\begin{proof}
As the same as \cite{LW19}, we let $\mathcal{A}_{\bar{\varphi}}(\varphi)=\Delta_\varphi\varphi+ \mathcal{L}_{V_{\bar{\varphi}}(\varphi)}\varphi$. Since $\tilde{\varphi}$ is closed, we have
$$\Delta_{\tilde{\varphi}}\tilde{\varphi}=dd^{\ast}_{\tilde{\varphi}}\tilde{\varphi}.$$
From $d\tilde{\varphi}=0$ and Cartan's formula, we see that
$$\mathcal{L}_{V_{\bar{\varphi}}(\tilde{\varphi})}\tilde{\varphi}=d(V_{\bar{\varphi}}(\tilde{\varphi})\lrcorner\tilde{\varphi}).$$

Thus, we get
$$\mathcal{A}_{\bar{\varphi}}(\tilde{\varphi})=d\left(d_{\widetilde\varphi}^{\ast}\widetilde\varphi+V_{\bar{\varphi}}(\widetilde\varphi)\lrcorner\widetilde\varphi\right).$$
The codifferential operator tells us that there exists a linear map
$$B_{1}:T_{x}^{\ast}M\otimes \Lambda^{3}T_{x}^{\ast}M\rightarrow \Lambda^{2}T_{x}^{\ast}M,\quad \text{for all}\ x\in M$$
such that 
$$d_{\widetilde\varphi}^{\ast}\widetilde\varphi=B_{1}(x,\theta)\bar{\nabla}\theta,$$
where we use $\bar{\nabla}\tilde{\varphi}=\bar{\nabla}\theta$.

From (4.9)-(4.12) in \cite{LW19}, we know that there also exists a linear map $B_{2}(x,\theta)$ such that
$$V_{\bar{\varphi}}(\widetilde\varphi)\lrcorner\widetilde\varphi=B_{2}(x,\theta)\bar{\nabla}\theta.$$
By setting $B(x,\theta)=B_{1}(x,\theta)+B_{2}(x,\theta)$, we get
$$d_{\widetilde\varphi}^{\ast}\widetilde\varphi+V_{\bar{\varphi}}(\widetilde\varphi)\lrcorner\widetilde\varphi=B(x,\theta)\bar{\nabla}\theta,$$
which implies
\begin{align}\label{A closed}
    \mathcal{A}_{\bar{\varphi}}(\tilde{\varphi})=d\left(B(x,\theta)\bar{\nabla}\theta\right).
\end{align}

Now, we let $\varphi_{s}=\bar{\varphi}+s\theta$, then $d\theta=d(\bar{\varphi}-\tilde{\varphi})=0$ yields $d\varphi_{s}=0$. Bryant and Xu \cite{BX11} show that the linearisation of $\mathcal{A}_{\bar{\varphi}}(\varphi_{s})$ is
$$\left.\frac{d}{ds}\right|_{s=0}\mathcal{A}_{\bar{\varphi}}(\varphi_{s})=-\Delta_{\bar{\varphi}}\theta,$$
which is also (4.5) in \cite{LW19}. From \eqref{A closed}, we also get
\begin{align}\label{Delta}
    \left.\frac{d}{ds}\right|_{s=0}\mathcal{A}_{\bar{\varphi}}(\varphi_{s})=d\left(B(x,0)\bar{\nabla}\theta\right)
\end{align}
Together with \eqref{A closed} and \eqref{Delta}, we show that
\begin{align}
    \frac{\partial}{\partial t}\theta(t)=-\Delta_{\bar{\varphi}}\theta(t)+dF(\bar{\varphi},\tilde{\varphi}(t),\theta(t),\bar{\nabla}\theta(t)),
\end{align}
where $F(\bar{\varphi},\tilde{\varphi}(t),\theta(t),\bar{\nabla}\theta(t))=\left(B(x,\theta)-B(x,0)\right)\bar{\nabla}\theta$.

If $\|\theta(t)\|_{C^{0}_{\bar{g}}}<\delta$ for a small $\delta$, the open condition of the $G_{2}$-structures shows that $B$ is smooth and $D_{z}B(x,z)$ is uniformly bounded in the $\delta$-neighborhood of $\bar{\varphi}$. Thus, we get
$$B(x,\theta)-B(x,0)=\int_{0}^{1}D_{z}B(x,s\theta)[\theta]ds.$$
It follows that $|F|_{\bar{g}}\leq C|\theta(t)|_{\bar{g}}|\bar{\nabla}\theta(t)|_{\bar{g}}$.
\end{proof}

For the Hodge Laplacian operator $\Delta$, if $\Delta\eta_{i}=\lambda_{i}\eta_{i}$ for some form $\eta_{i}$, we call $\lambda_{i}$ an eigenvalue of $\Delta$, and $\eta_{i}$ is the eigenform with respect to $\lambda_{i}$. For any smooth functions $a_{i}$, if $\eta=\sum_{i}a_{i}\eta_{i}$, then we write
$$e^{-t\Delta}\eta=\sum_{i}e^{-t\lambda_{i}}a_{i}\eta_{i}.$$
By semigroup theory, we obtain the following lemma.
\begin{lemma}\label{lem3.2}
For all time $t$ for which the solution $\tilde{\varphi}(t)$ of \eqref{theta} exists, we have
\begin{align}
    \theta(t)=e^{-t\Delta_{\bar{\varphi}}}\theta(0)+\int_{0}^{t}de^{-(t-s)\Delta_{\bar{\varphi}}}F(s)ds.
\end{align}
\end{lemma}
\begin{proof}
Fix $t>0$, we let $0\leq s\leq t$, differentiating with respect to 
$s$ on the heat operator act on $\theta$, we obtain
\begin{align}\label{heat operator}
    \frac{d}{ds}\left(e^{-(t-s)\Delta_{\bar{\varphi}}}\theta(s)\right)&=e^{-(t-s)\Delta_{\bar{\varphi}}}(\Delta_{\bar{\varphi}}\theta(s)+\partial_{s}\theta(s))\\
    &=e^{-(t-s)\Delta_{\bar{\varphi}}}dF(s)\notag.
\end{align}
By integrating \eqref{heat operator} from $0$ to $t$, we get
$$\theta(t)-e^{-t\Delta_{\bar{\varphi}}}\theta(0)=\int_{0}^{t}e^{-(t-s)\Delta_{\bar{\varphi}}}dF(s)\ ds,$$
which implies
\begin{align}
    \theta(t)&=e^{-t\Delta_{\bar{\varphi}}}\theta(0)+\int_{0}^{t}e^{-(t-s)\Delta_{\bar{\varphi}}}dF(s)\ ds\notag\\
    &=e^{-t\Delta_{\bar{\varphi}}}\theta(0)+\int_{0}^{t}de^{-(t-s)\Delta_{\bar{\varphi}}}F(s)\ ds\notag.
\end{align}
The second equality uses Lemma 4.17 in \cite{Ros97}.
\end{proof}

\begin{lemma}\label{lem3.3}
    Let $(M,\bar{g})$ be a compact $7$-dimensional Riemannian manifold, $\bar{g}=g_{\bar{\varphi}}$ is the Riemannian metric induced by the torsion-free $G_{2}$-structure $\bar{\varphi}$. Suppose that $\eta$ is a $p$-form satisfying
    \begin{align}
        \frac{\partial}{\partial t}\eta+\Delta_{\bar{\varphi}}\eta=0,\quad \quad \eta(0)=\eta_{0},\quad\quad t\in[0,1].
    \end{align}
    Let $0<\alpha<1$ and $k\in\mathbb{N}$, then there exists a uniform constant $C$ depend only on $M,\bar{\varphi},k,\alpha$, such that for all $t\in(0,1]$ 
    \begin{align}\label{3.8}
        \|\bar{\nabla}^{k}\eta(t)\|_{C^{0}}\leq Ct^{-\frac{k}{2}}\|\eta_{0}\|_{C^{0}},&\quad \quad [\bar{\nabla}^{k}\eta(t)]_{\alpha}\leq Ct^{-\frac{k+\alpha}{2}}\|\eta_{0}\|_{C^{0}};\\\label{3.9}
        \|\bar{\nabla}^{k}d\eta(t)\|_{C^{0}}\leq Ct^{-\frac{k+1}{2}}\|\eta_{0}\|_{C^{0}},& \quad \quad[\bar{\nabla}^{k}d\eta(t)]_{\alpha}\leq Ct^{-\frac{k+1+\alpha}{2}}\|\eta_{0}\|_{C^{0}}.
    \end{align}
    In addition, for $k\geq 1$
    \begin{align}
        \|\bar{\nabla}^{k}d\eta(t)\|_{C^{0}}&\leq Ct^{-1+\frac{\alpha}{2}}\|\eta_{0}\|_{C^{k-1,\alpha}},\notag\\
        \|\bar{\nabla}^{k+1}d\eta(t)\|_{C^{0}}&\leq Ct^{-\frac{3}{2}+\frac{\alpha}{2}}\|\eta_{0}\|_{C^{k-1,\alpha}}.\notag
    \end{align}
\end{lemma}
\begin{proof}
    Recall the Weitzenb\"ock formula:
    $$\Delta_{\bar{\varphi}}\eta=-\triangle_{\bar{g}}\eta+\mathcal{R}(\eta),$$
    where $\mathcal{R}$ is some combination of curvature operators with respect to $\bar{g}$ and $\triangle_{\bar{g}}$ is the connection Laplacian. It follows that
    \begin{align}
    \left(\frac{\partial}{\partial t}-\triangle_{\bar{g}}\right)\eta=-\mathcal{R}(\eta).\notag
    \end{align}
    Then, by calculation, we get
    \begin{align}\label{0order}
        \left(\partial_{t}-\triangle_{\bar{g}}\right)|\eta|^{2}&=2\langle(\partial_t-\triangle_{\bar{g}})\eta,\eta\rangle-2|\bar\nabla\eta|^2\\
        &=-2\langle\mathcal R(\eta),\eta\rangle-2|\bar\nabla\eta|^2\notag\\
        &\leq -2|\bar\nabla\eta|^2+C|\eta|^{2}.\notag
    \end{align}
    For derivatives of order $k$, we have
    \begin{align}
        (\partial_{t}-\triangle_{\bar{g}})\bar{\nabla}^{k}\eta&=\bar{\nabla}^{k}\left((\partial_{t}-\triangle_{\bar{g}})\eta\right)+\sum_{i=0}^{k}B_{k,i}(\mathcal{R})\bar{\nabla}^{i}\eta\notag\\
        &=\sum_{i=0}^{k}B_{k,i}(\mathcal{R})\bar{\nabla}^{i}\eta.\notag
    \end{align}
    where $B_{k,i}(\mathcal{R})$ is a smooth uniformly polynomial consisting of the background curvature and its derivatives up to finite order. Similar to \eqref{0order}, we have
    \begin{align}\label{korder}
         (\partial_{t}-\triangle_{\bar{g}})|\bar{\nabla}^{k}\eta|^{2}\leq -2|\bar{\nabla}^{k+1}\eta|^{2}+C\sum_{i=0}^{k}|\bar{\nabla}^{i}\eta|^{2}.
    \end{align}
    For $0\leq j\leq k$, we define $$E_{k,j}:=\sum_{i=0}^{j}|\bar{\nabla}^{i}\eta(t)|^{2}+\sum_{i=j+1}^{k}\frac{t^{i-j}}{(i-j)!}|\bar{\nabla}^{i}\eta|^{2}.$$
    Then by \eqref{0order} and \eqref{korder}, we have
    \begin{align}
        (\partial_{t}-\triangle_{\bar{g}})E_{k,j}&=\sum_{i=0}^{j}(\partial_{t}-\triangle_{\bar{g}})|\nabla^{i}\eta(t)|^{2}+\sum_{i=j+1}^{k}\frac{t^{i-j}}{(i-j)!}(\partial_{t}-\triangle_{\bar{g}})|\bar{\nabla}^{i}\eta|^{2}\notag\\
        &\quad+\sum_{i=j+1}^{k}\frac{t^{i-j-1}}{(i-j-1)!}|\bar{\nabla}^{i}\eta|^{2}\notag\\
        &\leq -2\sum_{i=0}^{j}|\nabla^{i+1}\eta(t)|^{2}-2\sum_{i=j+1}^{k}\frac{t^{i-j}}{(i-j)!}|\bar{\nabla}^{i+1}\eta|^{2}+C\sum_{i=0}^{j}|\nabla^{i}\eta(t)|^{2}\notag\\
        &\quad+C\sum_{i=j+1}^{k}\frac{t^{i-j}}{(i-j)!}|\bar{\nabla}^{i}\eta|^{2}+\sum_{i=j+1}^{k}\frac{t^{i-j-1}}{(i-j-1)!}|\bar{\nabla}^{i}\eta|^{2}\notag\\
        &\leq CE_{k,j}.\notag
    \end{align}
    By setting $\hat{E}_{k,j}=e^{-Ct}E_{k,j}$, we get
    $$ (\partial_{t}-\triangle_{\bar{g}})\hat{E}_{k,j}\leq 0.$$
    Then the maximum principle tells us that
    \begin{align}
        E_{k,j}(x,t)\leq Ce^{C\times 0}E_{k,j}(x,0)=C\sum_{i=0}^{j}|\bar{\nabla}^{i}\eta(0)|^{2},\quad \quad x\in M.\notag
    \end{align}
    which implies
    \begin{align}\label{E_k,j}
        \|\bar{\nabla}^{k}\eta(t)\|_{C^{0}}\leq Ct^{-\frac{k-j}{2}}\|\eta_{0}\|_{C^{j}}.
    \end{align}
    Moreover, if we let $j=0$, then $\eqref{E_k,j}$ implies $ \|\bar{\nabla}^{k}\eta(t)\|_{C^{0}}\leq Ct^{-\frac{k}{2}}\|\eta_{0}\|_{C^{0}}$.

    Let $0<\rho<{\rm inj}(M,\bar{g})$, for any $x,y\in M$ satisfies $d_{\bar{g}}(x,y)< \rho$, we connect $x$ to $y$ by a geodesic $\gamma$. For any $p$-tensor $S$, we have
    \begin{align}
        |S(x)-P_{\gamma}S(y)|\leq 2\|S\|_{C^{0}}\notag
    \end{align}
    and 
    \begin{align}
        |S(x)-P_{\gamma}S(y)|\leq \|\bar{\nabla}S\|_{C^{0}}\cdot d_{\bar{g}}(x,y)\notag.
    \end{align}
    It follows that
    \begin{align}
        |S(x)-P_{\gamma}S(y)|\leq \min\left\{ 2\|S\|_{C^{0}},\|\bar{\nabla}S\|_{C^{0}}\cdot d_{\bar{g}}(x,y)\right\}.\notag
    \end{align}
    Using interpolation inequality $\min\{a,b\}\leq a^{1-\alpha}b^{\alpha}$ with $\alpha\in(0,1)$, we get
    \begin{align}\label{inter}
        [S]_{\alpha}\leq C\frac{\|S\|_{C^{0}}^{1-\alpha}\|\bar{\nabla}S\|_{C^{0}}^{\alpha} d_{\bar{g}}(x,y)^{\alpha}}{d_{\bar{g}}(x,y)^{\alpha}}= C\|S\|_{C^{0}}^{1-\alpha}\|\bar{\nabla}S\|_{C^{0}}^{\alpha}.
    \end{align}
    Substituting $S=\bar{\nabla}^{k}\eta$ into \eqref{inter}, we obtain
    \begin{align}
        [\bar{\nabla}^{k}\eta]_{\alpha}&\leq \|\bar{\nabla}^{k}\eta\|_{C^{0}}^{1-\alpha}\|\bar{\nabla}^{k+1}\eta\|_{C^{0}}^{\alpha}\notag\\
        &\leq C\left(t^{-\frac{k}{2}}\|\eta_{0}\|_{C^{0}}\right)^{1-\alpha}\cdot \left(t^{-\frac{k+1}{2}}\|\eta_{0}\|_{C^{0}}\right)^{\alpha}\notag\\
        &\leq Ct^{-\frac{k+\alpha}{2}}\|\eta_{0}\|_{C^{0}}.\notag
    \end{align}
    Thus we prove \eqref{3.8}. \eqref{3.9} can be obtained by $|d\eta|\leq C|\bar{\nabla}\eta|$ and \eqref{3.8}.

Choose finitely many global smooth vector fields $X_\nu$ which
agree with coordinate vector fields on neighborhoods of the
closures of smaller coordinate charts covering $M$.
Such fields are obtained using fixed smooth cutoffs.
Let $\Phi_h^\nu$ be their flows. Fix $h_0\in(0,1]$ sufficiently
small that these flows act as coordinate translations on the
corresponding smaller charts for $0<h\leq h_0$. Set
$$
T_h^\nu=(\Phi_h^\nu)^{\ast},
\qquad \delta_h^\nu\eta=T_h^\nu\eta-\eta,
\qquad \Delta_h^\nu=T_h^\nu\Delta_{\bar\varphi}(T_h^\nu)^{-1}.
$$
The initial H\"older bound gives
$$
\|\delta_h^\nu\eta_0\|_{C^{j-1}}\leq Ch^\alpha\|\eta_0\|_{C^{j-1,\alpha}}.
$$
Taking the difference of the heat equations yields
$$\partial_t(\delta_h^\nu\eta)+\Delta_h^\nu(\delta_h^\nu\eta)=
(\Delta_{\bar\varphi}-\Delta_h^\nu)\eta.$$
The coefficients of $\Delta_h^\nu$ depend smoothly on $h$, and
$\Delta_0^\nu=\Delta_{\bar\varphi}$. Hence
$$
\|(\Delta_{\bar\varphi}-\Delta_h^\nu)\eta\|_{C^{j-1}}
\leq Ch\|\eta\|_{C^{j+1}}.
$$
Moreover,
$$
e^{-t\Delta_h^\nu}
=
T_h^\nu e^{-t\Delta_{\bar\varphi}}(T_h^\nu)^{-1},
$$
All these constants are independent of $h$ and $\nu$. Similar to Lemma \ref{lem3.2}, Duhamel's formula therefore gives
$$\max_\nu\|\delta_h^\nu\eta(t)\|_{C^{j-1}}\leq
Ch^\alpha\|\eta_0\|_{C^{j-1,\alpha}}+Ch\int_0^t\|\eta(s)\|_{C^{j+1}}\,ds.
$$

For a scalar function in a coordinate chart, the fundamental
theorem of calculus gives
$$\partial_l f(x)=\frac{f(x+he_l)-f(x)}{h}-\frac1h\int_0^h(h-r)
\partial_l^2f(x+re_l)\,dr,$$
where $e_l$ is the corresponding coordinate unit vector.
Apply this identity to the coordinate derivatives of order
$j-1$ of each coefficient of $\eta(t)$. Since the smaller
charts cover $M$, comparison of coordinate and covariant
derivatives gives
\begin{align}
\|\eta(t)\|_{C^j}\leq Ch^{-1}\max_\nu
\|\delta_h^\nu\eta(t)\|_{C^{j-1}}+Ch\|\eta(t)\|_{C^{j+1}}
+C\|\eta(t)\|_{C^{j-1}}.\notag
\end{align}

The last term accounts for the lower-order terms in this
comparison. Substituting the difference estimate and using
$$\|\eta(t)\|_{C^{j-1}}\leq C\|\eta_0\|_{C^{j-1}},
\qquad h^{\alpha-1}\geq1,$$
we obtain
\begin{align}
\|\eta(t)\|_{C^j}\leq Ch^{\alpha-1}\|\eta_0\|_{C^{j-1,\alpha}}+C\int_0^t\|\eta(s)\|_{C^{j+1}}\,ds+Ch\|\eta(t)\|_{C^{j+1}}.
\end{align}

For $0<T\leq1$, define
$$\beta(T)=\sup_{0<t\leq T}t^{(1-\alpha)/2}\|\eta(t)\|_{C^j}.$$
Since $\eta_0$ is smooth, $\beta(T)$ is finite.
By the semigroup property and the integer-order estimate,
$$\|\eta(s)\|_{C^{j+1}}
\leq Cs^{-1/2}\|\eta(s/2)\|_{C^j}
\leq C\beta(T)s^{-1+\alpha/2},
\qquad 0<s\leq T.
$$
Since $\alpha>0$, it follows that
$$
\int_0^t\|\eta(s)\|_{C^{j+1}}\,ds
\leq C\beta(T)t^{\alpha/2}.
$$
Take $h=\varepsilon\sqrt t$, with
$0<\varepsilon\leq h_0$, in the preceding difference-quotient
estimate. Multiplying by $t^{(1-\alpha)/2}$ and taking the
supremum over $0<t\leq T$ gives
$$\beta(T)\leq C_0\varepsilon^{\alpha-1}\|\eta_0\|_{C^{j-1,\alpha}}
+C_1(\sqrt T+\varepsilon)\beta(T),$$
where $C_0$ and $C_1$ are independent of $\varepsilon$, $T$,
and $\eta_0$. Choose $\varepsilon$ so that
$C_1\varepsilon\leq1/4$, and then choose $T_0\in(0,1]$ so that
$C_1\sqrt{T_0}\leq1/4$. Absorbing the last term yields
$$\beta(T_0)\leq 2C_0\varepsilon^{\alpha-1}\|\eta_0\|_{C^{j-1,\alpha}}.
$$
For $T_0\leq t\leq1$, the integer-order estimate gives
$$\|\eta(t)\|_{C^j}\leq Ct^{-1/2}\|\eta_0\|_{C^{j-1}}
\leq CT_0^{-\alpha/2}t^{-1/2+\alpha/2}\|\eta_0\|_{C^{j-1,\alpha}}.$$
Consequently, after enlarging $C$,
$$\|\eta(t)\|_{C^j}\leq Ct^{-1/2+\alpha/2} \|\eta_0\|_{C^{j-1,\alpha}},
\qquad 0<t\leq1.$$

Finally, exterior differentiation is the alternation of
$\bar\nabla$. Applying the integer-order estimate on the
interval $[t/2,t]$, we obtain
\begin{align}
\|\bar\nabla^j d\eta(t)\|_{C^0}
&\leq Ct^{-1/2}\|\eta(t/2)\|_{C^j}\leq
Ct^{-1+\alpha/2}\|\eta_0\|_{C^{j-1,\alpha}},\\
\|\bar\nabla^{j+1}d\eta(t)\|_{C^0}
&\leq Ct^{-1}\|\eta(t/2)\|_{C^j}\leq
Ct^{-3/2+\alpha/2}\|\eta_0\|_{C^{j-1,\alpha}},
\qquad 
\end{align}
for $0< t \leq1$. All constants depending only on the fixed background, $j$ and $\alpha$.
\end{proof}

\section{Main Theorem}\label{section4}
In this section, we derive an estimate for $\theta(t)$, and based on this estimate, we prove the main estimate.

Firstly, we define 
\begin{align}\label{K}
    K_{n}(\tau):=\sup_{0<t\leq \tau}\left\{\|\theta(t)\|_{C^{0}_{\bar{g}}}+\sum_{i=1}^{n}t^{\frac{i}{2}}\|\bar{\nabla}^{i}\theta(t)\|_{C^{0}_{\bar{g}}}+\sum_{i=1}^{n}t^{\frac{i+\alpha}{2}}[\bar{\nabla}^{i}\theta(t)]_{\alpha}\right\}.
\end{align}
Then we get the following Lemma:
\begin{lemma}\label{lem4.1}
    For all time $t$ for which the solution $\tilde{\varphi}(t)$ of \eqref{theta} exists. There exists $\delta^{\ast}=\delta^{\ast}(\bar{\varphi})>0$ such that if $\delta\in(0,\delta^{\ast})$ and $\|\theta(0)\|_{C^{0}_{\bar{g}}}<\delta$, then for $0<\alpha<1$ and $K_{n}(\tau)\leq 1$,
    \begin{align}\label{Kinequality}
        K_{n}(\tau)\leq C_{0}\delta+C_{1}K_{n}(\tau)^{2},\quad \quad 0<\tau\leq 1.
    \end{align}
    where $C_{0}$ and $C_{1}$ are uniformly  constants depending only on $M,\bar{\varphi},n,\alpha$.
\end{lemma}
\begin{proof}
    Fix $\tau$, by the definition of $K_{n}=K_{n}(\tau)$, we know that
    \begin{align}\label{4.3}
        \|\theta(t)\|_{C^{0}}\leq K_{n},\quad \quad \|\bar{\nabla}^{i}\theta(t)\|_{C^{0}}\leq K_{n}t^{-\frac{i}{2}},\quad\quad [\bar{\nabla}^{i}\theta(t)]_{\alpha}\leq K_{n}t^{-\frac{i+\alpha}{2}}
    \end{align}
    for all $1\leq i\leq n$ and $t\in(0,\tau]$. By the interpolation inequality in Lemma \ref{lem3.3}, we have
    \begin{align}
        [\theta(t)]_{\alpha}\leq C \|\theta(t)\|_{C^{0}}^{1-\alpha}\|\bar{\nabla}\theta(t)\|_{C^{0}}^{\alpha}\leq CK_{n}t^{-\frac{\alpha}{2}}.
    \end{align}

    Now we prove \eqref{Kinequality}. From the definition of $K_{n}$, we should estimate three parts in \eqref{K}. From Lemma \ref{lem3.2}, we should estimate $F(t)$ and its derivative firstly. By Lemma \ref{lem3.1} and \eqref{4.3}, we have
    \begin{align}
        \|F(t)\|_{C^{0}}\leq C\|\theta(t)\|_{C^{0}}\|\bar{\nabla}\theta(t)\|_{C^{0}}\leq CK_{n}^{2}t^{-\frac{1}{2}}.\notag
    \end{align}
    Lemma \ref{lem3.1} tells us that $F(t)=f(\bar{\varphi},\tilde{\varphi}(t))\theta(t)\bar{\nabla}\theta(t)$, combine with $\tilde{\varphi}=\bar{\varphi}-\theta$ and $\bar{\nabla}\tilde{\varphi}=\bar{\nabla}\theta$, we know that each term in the expansion of $\bar{\nabla}^{i}F(t)$ can be expressed as 
    $$G(\bar{\varphi},\theta(t))\ast\bar{\nabla}^{a_{1}}\theta\ast\cdots\ast \bar{\nabla}^{a_{q}}\theta$$
    with $1\leq i\leq n-1$, $q\geq 2, a_{j}\geq 0$ and $a_{1}+\cdots+a_{q}= i+1$. Since $G(\bar{\varphi},\theta(t))$ is a smooth function, we obtain
    \begin{align}
        \|\bar{\nabla}^{i}F(t)\|_{C^{0}}&\leq C\|\bar{\nabla}^{a_{1}}\theta(t)\|_{C^{0}}\cdots \|\bar{\nabla}^{a_{q}}\theta(t)\|_{C^{0}}\\
        &\leq CK_{n}^{q}t^{-\frac{a_{1}+\cdots+a_{q}}{2}}\notag\\
        &\leq CK_{n}^{2}t^{-\frac{i+1}{2}}\notag
    \end{align}
    for $1\leq i\leq n-1$. The last inequality uses $K_{n}\leq 1$.

    For any $C^{\alpha}$ function $u_{1},\cdots,u_{q}$, recall the product inequality for H\"older norms, we get
$$[u_{1}\cdots u_{q}]_{\alpha}\leq \sum_{i=1}^{q}[u_{i}]_{\alpha}\prod_{i\neq j}\|u_{j}\|_{C^{0}}.$$
Then the H\"older norm of $\bar{\nabla}^{i}F(t)$ satisfies
   \begin{align}\label{i-1}
    [\bar{\nabla}^{i}F(t)]_{\alpha}&\leq C[\bar{\nabla}^{a_{1}}\theta\ast\cdots\ast \bar{\nabla}^{a_{q}}\theta]_{\alpha}\\
    &\leq C\sum_{l=1}^{q}[\bar{\nabla}^{a_{l}}\theta]_{\alpha}\prod_{l\neq j}\|\bar{\nabla}^{a_{j}}\theta\|_{C^{0}}\notag\\
    &\leq CK_{n}^{q}t^{-\frac{a_{1}+\cdots+a_{q}+\alpha}{2}}\notag\\
    &\leq  CK_{n}^{2}t^{-\frac{i+1+\alpha}{2}}\notag
   \end{align}
    for $1\leq i\leq n-1$. The last inequality also uses $K_{n}\leq 1$.

    Note that
    \begin{align}
        \left(\frac{\partial}{\partial t}+\Delta_{\bar{\varphi}}\right)(e^{-t\Delta_{\bar{\varphi}}}\theta(0))&=0,\quad \quad (e^{-t\Delta_{\bar{\varphi}}}\theta(0))(0)=\theta(0).\notag\\
        \left(\frac{\partial}{\partial (t-s)}+\Delta_{\bar{\varphi}}\right)(e^{-(t-s)\Delta_{\bar{\varphi}}}F(s))&=0,\quad \quad (e^{-(t-s)\Delta_{\bar{\varphi}}}F(s))|_{t=s}=F(s).\notag
    \end{align}
    for $t\in(0,\tau]$ and $0<s\leq t$. By Lemma \ref{lem3.3}, it follows that
    \begin{align}\label{4.6}
        \|\bar{\nabla}^{j}e^{-t\Delta_{\bar{\varphi}}}\theta(0)\|_{C^{0}}&\leq Ct^{-\frac{j}{2}}\|\theta(0)\|_{C^{0}}\leq Ct^{-\frac{j}{2}}\delta,\\ \label{4.7}
        [\bar{\nabla}^{j}e^{-t\Delta_{\bar{\varphi}}}\theta(0)]_{\alpha}&\leq Ct^{-\frac{j+\alpha}{2}}\|\theta(0)\|_{C^{0}}\leq Ct^{-\frac{j+\alpha}{2}}\delta,\\ \label{4.8}
        \|\bar{\nabla}^{j}de^{-(t-s)\Delta_{\bar{\varphi}}}F(s)\|_{C^{0}}&\leq C(t-s)^{-\frac{j+1}{2}}\|F(s)\|_{C^{0}},\\ \label{4.9}
        \|\bar{\nabla}^{j}de^{-(t-s)\Delta_{\bar{\varphi}}}F(s)\|_{C^{0}}&\leq C(t-s)^{-1+\frac{\alpha}{2}}\|F(s)\|_{C^{j-1,\alpha}},\\ \label{4.10}
        \|\bar{\nabla}^{j+1}de^{-(t-s)\Delta_{\bar{\varphi}}}F(s)\|_{C^{0}}&\leq C(t-s)^{-\frac{3}{2}+\frac{\alpha}{2}}\|F(s)\|_{C^{j-1,\alpha}}.
    \end{align}

    Now, we can estimate the part $\|\theta(t)\|_{C^{0}}$. From \eqref{4.6}, we know that
    $$\|e^{-t\Delta_{\bar{\varphi}}}\theta(0)\|_{C^{0}}\leq C\delta$$
    and \eqref{4.8} tells us that
    \begin{align}
        \left\|\int_{0}^{t}de^{-(t-s)\Delta_{\bar{\varphi}}}F(s)ds\right\|_{C^{0}}&\leq \int_{0}^{t}\left\|de^{-(t-s)\Delta_{\bar{\varphi}}}F(s)\right\|_{C^{0}}ds\notag\\
        &\leq C\int_{0}^{t}(t-s)^{-\frac{1}{2}}\|F(s)\|_{C^{0}}ds\notag\\
        &\leq CK_{n}^{2}\int_{0}^{t}(t-s)^{-\frac{1}{2}}s^{-\frac{1}{2}}ds\notag
    \end{align}
    By setting $s=t\sin^{2}u$, then $0\leq u\leq \frac{\pi}{2}$, we get
    $$\int_{0}^{t}(t-s)^{-\frac{1}{2}}s^{-\frac{1}{2}}ds=\pi.$$
    Thus we get 
    \begin{align}\label{C0}
        \|\theta(t)\|_{C^{0}}\leq C\delta+CK_{n}^{2}.
    \end{align}

    In the following, we estimate $t^{\frac{i}{2}}\|\bar{\nabla}^{i}\theta(t)\|_{C^{0}_{\bar{g}}}$ for $1\leq i\leq n$. \eqref{4.6} yields
    \begin{align}
        t^{\frac{i}{2}}\|\bar{\nabla}^{i}e^{-t\Delta_{\bar{\varphi}}}\theta(0)\|_{C^{0}}\leq C\delta,\quad \quad 1\leq i\leq n.\notag
    \end{align}
    For the term contained $F$, we discuss it in $2$ parts. When $0<s\leq\frac{t}{2}$, \eqref{4.8} yields
    \begin{align}
        \left\|\bar{\nabla}^{i}\int_{0}^{\frac{t}{2}}de^{-(t-s)\Delta_{\bar{\varphi}}}F(s)ds\right\|_{C^{0}}&\leq \int_{0}^{\frac{t}{2}}\left\|\bar{\nabla}^{i}de^{-(t-s)\Delta_{\bar{\varphi}}}F(s)\right\|_{C^{0}}ds\notag\\
        &\leq C\int_{0}^{\frac{t}{2}}(t-s)^{-\frac{i+1}{2}}\|F(s)\|_{C^{0}}ds\notag\\
        &\leq CK_{n}^{2}\int_{0}^{\frac{t}{2}}(t-s)^{-\frac{i+1}{2}}s^{-\frac{1}{2}}ds\notag\\
        &\leq CK_{n}^{2}\left(\frac{t}{2}\right)^{-\frac{i+1}{2}}\int_{0}^{\frac{t}{2}}s^{-\frac{1}{2}}ds\notag\\
        &\leq CK_{n}^{2}t^{-\frac{i}{2}}\notag.
    \end{align}
    The fourth inequality uses $t-s\geq \frac{t}{2}$.

    When $\frac{t}{2}<s\leq t$, by \eqref{i-1} and \eqref{4.9}, we obtain
    \begin{align}
        \left\|\bar{\nabla}^{i}\int_{\frac{t}{2}}^{t}de^{-(t-s)\Delta_{\bar{\varphi}}}F(s)ds\right\|_{C^{0}}&\leq \int_{\frac{t}{2}}^{t}\left\|\bar{\nabla}^{i}de^{-(t-s)\Delta_{\bar{\varphi}}}F(s)\right\|_{C^{0}}ds\notag\\
        &\leq C\int_{\frac{t}{2}}^{t}(t-s)^{-1+\frac{\alpha}{2}}\|F(s)\|_{C^{i-1,\alpha}}ds\notag\\
        &\leq CK_{n}^{2}\int_{\frac{t}{2}}^{t}(t-s)^{-1+\frac{\alpha}{2}}s^{-\frac{i+\alpha}{2}}ds\notag\\
        &\leq CK_{n}^{2}t^{-\frac{i+\alpha}{2}}\int_{\frac{t}{2}}^{t}(t-s)^{-1+\frac{\alpha}{2}}ds\notag\\
        &\leq CK_{n}^{2}t^{-\frac{i+\alpha}{2}}t^{\frac{\alpha}{2}}\notag\\
        &=CK_{n}^{2}t^{-\frac{i}{2}}\notag.
    \end{align}
    Together with the above, we get
    \begin{align}
        t^{\frac{i}{2}}\|\bar{\nabla}^{i}\theta(t)\|_{C^{0}_{\bar{g}}}\leq C\delta+CK_{n}^{2}\notag,
    \end{align}
    which implies
    \begin{align}\label{Ci}
        \sum_{i=1}^{n}t^{\frac{i}{2}}\|\bar{\nabla}^{i}\theta(t)\|_{C^{0}_{\bar{g}}}\leq C\delta+CK_{n}^{2}.
    \end{align}

    Lastly, we estimate the term $t^{\frac{i+\alpha}{2}}[\bar{\nabla}^{i}\theta(t)]_{\alpha}$ for $1\leq i\leq n$. From \eqref{4.7}, it follows that
    \begin{align}\label{alpha1}
        t^{\frac{i+\alpha}{2}}[\bar{\nabla}^{i}e^{-t\Delta_{\bar{\varphi}}}\theta(0)]_{\alpha}\leq C\delta.
    \end{align}
    Now, we estimate the H\"older term contained $F$. We also estimate it into two parts. When $0<s\leq \frac{t}{2}$, by \eqref{4.8} and the interpolation inequality in Lemma \ref{lem3.3}, we get
    \begin{align}\label{alpha2}
        &\quad\left[\bar{\nabla}^{i}\int_{0}^{\frac{t}{2}}de^{-(t-s)\Delta_{\bar{\varphi}}}F(s)ds\right]_{\alpha}\\
        &\leq C\left\|\bar{\nabla}^{i}\int_{0}^{\frac{t}{2}}de^{-(t-s)\Delta_{\bar{\varphi}}}F(s)ds\right\|^{1-\alpha}\cdot\left\|\bar{\nabla}^{i+1}\int_{0}^{\frac{t}{2}}de^{-(t-s)\Delta_{\bar{\varphi}}}F(s)ds\right\|^{\alpha}\notag\\
        &\leq C\left(K_{n}^{2}t^{-\frac{i}{2}}\right)^{1-\alpha}\cdot \left(K_{n}^{2}t^{-\frac{i+1}{2}}\right)^{\alpha}\notag\\
        &=CK_{n}^{2}t^{-\frac{i+\alpha}{2}}\notag.
    \end{align}
    When $\frac{t}{2}<s\leq t$, we let $u=t-s$, then $0<u\leq \frac{t}{2}$. Let $0<\rho<{\rm inj}(M,\bar{g})$, for any $x,y\in M$ satisfies $d_{\bar{g}}(x,y)< \rho$, we connect $x$ to $y$ by a geodesic $\gamma$. By \eqref{4.10}, we have
    \begin{align}
        &\quad\left|\left(\bar{\nabla}^{i}de^{-u\Delta_{\bar{\varphi}}}F(t-u)\right)(x)-P_{\gamma}\left(\bar{\nabla}^{i}de^{-u\Delta_{\bar{\varphi}}}F(t-u)\right)(y)\right|\notag\\
        &\leq C\left\|\bar{\nabla}^{i}de^{-u\Delta_{\bar{\varphi}}}F(t-u)\right\|_{C^{0}}\notag\\
        &\leq Cu^{-1+\frac{\alpha}{2}}\left\|F(t-u)\right\|_{C^{i-1,\alpha}}\notag.
    \end{align}
    Set $l=d_{\bar g}(x,y)$, by \eqref{4.9}, it follows that
    \begin{align}
        &\quad\left|\left(\bar{\nabla}^{i}de^{-u\Delta_{\bar{\varphi}}}F(t-u)\right)(x)-P_{\gamma}\left(\bar{\nabla}^{i}de^{-u\Delta_{\bar{\varphi}}}F(t-u)\right)(y)\right|\notag\\
        &\leq \left\|\bar{\nabla}^{i+1}de^{-u\Delta_{\bar{\varphi}}}F(t-u)\right\|_{C^{0}}\cdot l \notag\\
        &\leq Cu^{-\frac{3}{2}+\frac{\alpha}{2}}\left\|F(t-u)\right\|_{C^{i-1,\alpha}}\cdot l \notag.
    \end{align}
    which means
    \begin{align}
        &\quad\left|\left(\bar{\nabla}^{i}de^{-u\Delta_{\bar{\varphi}}}F(t-u)\right)(x)-P_{\gamma}\left(\bar{\nabla}^{i}de^{-u\Delta_{\bar{\varphi}}}F(t-u)\right)(y)\right|\notag\\
        &\leq C\left\|F(t-u)\right\|_{C^{i-1,\alpha}}\cdot\min\left\{u^{-1+\frac{\alpha}{2}},u^{-\frac{3}{2}+\frac{\alpha}{2}}\cdot l\right\}  \notag.
    \end{align}
    We know that $u=l^{2}$ if and only if $u^{-1+\frac{\alpha}{2}}=u^{-\frac{3}{2}+\frac{\alpha}{2}}\cdot l$, then we discuss this in parts.
    For $l^{2}<\frac{t}{2}$, by \eqref{i-1}, we have
    \begin{align}
         &\quad\left|\left(\bar{\nabla}^{i}\int_{0}^{\frac{t}{2}}de^{-u\Delta_{\bar{\varphi}}}F(t-u)du\right)(x)-P_{\gamma}\left(\bar{\nabla}^{i}\int_{0}^{\frac{t}{2}}de^{-u\Delta_{\bar{\varphi}}}F(t-u)du\right)(y)\right|\notag\\
        &\leq \int_{0}^{\frac{t}{2}}\left|\left(\bar{\nabla}^{i}de^{-u\Delta_{\bar{\varphi}}}F(t-u)\right)(x)-P_{\gamma}\left(\bar{\nabla}^{i}de^{-u\Delta_{\bar{\varphi}}}F(t-u)\right)(y)\right|du\notag\\
        &\leq C\int_{0}^{\frac{t}{2}}\left\|F(t-u)\right\|_{C^{i-1,\alpha}}\cdot\min\left\{u^{-1+\frac{\alpha}{2}},u^{-\frac{3}{2}+\frac{\alpha}{2}}\cdot l\right\}du\notag\\
        &\leq CK_{n}^{2}t^{-\frac{i+\alpha}{2}}\left(\int_{0}^{l^{2}}u^{-1+\frac{\alpha}{2}} du+\int_{l^{2}}^{\frac{t}{2}}u^{-\frac{3}{2}+\frac{\alpha}{2}}\cdot l\ du\right)\notag.
    \end{align}
    For $0<\alpha<1$, note that
    \begin{align}
        \int_{0}^{l^{2}}u^{-1+\frac{\alpha}{2}} du&=\frac{2}{\alpha}l^{\alpha}\notag,\\
        \int_{l^{2}}^{\frac{t}{2}}u^{-\frac{3}{2}+\frac{\alpha}{2}}\cdot l\ du&=\frac{2l}{1-\alpha}\left[\left(l^{2}\right)^{-\frac{1-\alpha}{2}}-\left(\frac{t}{2}\right)^{-\frac{1-\alpha}{2}}\right]\notag\\
        &\leq \frac{2l}{1-\alpha}\cdot \left(l^{2}\right)^{-\frac{1-\alpha}{2}}=\frac{2}{1-\alpha}l^{\alpha}\notag.
    \end{align}
    Thus, we get 
    \begin{align}
        &\quad\left|\left(\bar{\nabla}^{i}\int_{0}^{\frac{t}{2}}de^{-u\Delta_{\bar{\varphi}}}F(t-u)du\right)(x)-P_{\gamma}\left(\bar{\nabla}^{i}\int_{0}^{\frac{t}{2}}de^{-u\Delta_{\bar{\varphi}}}F(t-u)du\right)(y)\right|\notag\\
        &\leq CK_{n}^{2}\cdot t^{-\frac{i+\alpha}{2}}\cdot\left(\frac{2}{\alpha}+\frac{2}{1-\alpha}\right)l^{\alpha}\leq CK_{n}^{2}\cdot t^{-\frac{i+\alpha}{2}}l^{\alpha}\notag.
    \end{align}
    By the definition of H\"older norm, we get
    \begin{align}\label{alpha3}
        \left[\bar{\nabla}^{i}\int_{\frac{t}{2}}^{t}de^{-(t-s)\Delta_{\bar{\varphi}}}F(s)ds\right]_{\alpha}\leq CK_{n}^{2} t^{-\frac{i+\alpha}{2}},\quad l^{2}< \frac{t}{2}.
    \end{align}

    For $l^{2}\geq \frac{t}{2}$, by \eqref{i-1}, we obtain
    \begin{align}
         &\quad\left|\left(\bar{\nabla}^{i}\int_{0}^{\frac{t}{2}}de^{-u\Delta_{\bar{\varphi}}}F(t-u)du\right)(x)-P_{\gamma}\left(\bar{\nabla}^{i}\int_{0}^{\frac{t}{2}}de^{-u\Delta_{\bar{\varphi}}}F(t-u)du\right)(y)\right|\notag\\
        &\leq C\int_{0}^{\frac{t}{2}} \|F(t-u)\|_{C^{i-1,\alpha}}\cdot u^{-1+\frac{\alpha}{2}}du\notag\\
        &\leq CK_{n}^{2}t^{-\frac{i+\alpha}{2}}\int_{0}^{\frac{t}{2}}u^{-1+\frac{\alpha}{2}} du\notag\\
        &=CK_{n}^{2}t^{-\frac{i+\alpha}{2}}\cdot \frac{2}{\alpha}\left(\frac{t}{2}\right)^{\frac{\alpha}{2}}\notag\\
        &\leq CK_{n}^{2}t^{-\frac{i+\alpha}{2}} l^{\alpha}\notag
    \end{align}
    where the last inequality we use $l^{2}\geq \frac{t}{2}$. The above implies
    \begin{align}\label{alpha4}
        \left[\bar{\nabla}^{i}\int_{\frac{t}{2}}^{t}de^{-(t-s)\Delta_{\bar{\varphi}}}F(s)ds\right]_{\alpha}\leq CK_{n}^{2} t^{-\frac{i+\alpha}{2}},\quad l^{2}\geq \frac{t}{2}.
    \end{align}
    Then by \eqref{alpha1} to \eqref{alpha4}, we get
    \begin{align}
        \sum_{i=1}^{n}t^{\frac{i+\alpha}{2}}[\bar{\nabla}^{i}\theta(t)]_{\alpha}\leq C\delta+CK_{n}^{2}\notag.
    \end{align}
    Combing with \eqref{C0} and \eqref{Ci}, we prove \eqref{Kinequality}.
\end{proof}

\begin{theorem}
    For all time $t$ for which the solution $\tilde{\varphi}(t)$ of \eqref{theta} exists. If $\|\theta(0)\|_{C^{0}_{\bar{g}}}<\delta\leq \delta_{0}$ for a small $\delta_{0}>0$, then
    \begin{align}\label{main}
        K_{n}(\tau)\leq C\delta,\quad \quad 0<\tau< T_{\max},\quad \tau\leq 1,
    \end{align}
    and $T_{\max}> 1$,
    where $C$ is an uniformly  constant depend only on $\bar{\varphi}$, $\alpha$ and $n$.
\end{theorem}
\begin{proof}
    We define
    $$q_{n}(t):=\|\theta(t)\|_{C^{0}_{\bar{g}}}+\sum_{i=1}^{n}t^{\frac{i}{2}}\|\bar{\nabla}^{i}\theta(t)\|_{C^{0}_{\bar{g}}}+\sum_{i=1}^{n}t^{\frac{i+\alpha}{2}}[\bar{\nabla}^{i}\theta(t)]_{\alpha},$$
    then $K_{n}=\displaystyle\sup_{0<t\leq\tau}q_{n}(t)$. Since $\theta(t)$ is smooth, there exists a constant $A_{0}$ such that for a small time $t_{0}\in[0,T_{\max})$, 
    $$\sup_{0\leq t\leq t_{0}}\left(\sum_{i=1}^{n}\|\bar{\nabla}^{i}\theta(t)\|_{C^{0}_{\bar{g}}}+\sum_{i=1}^{n}[\bar{\nabla}^{i}\theta(t)]_{\alpha}\right)\leq A_{0}.$$
    Thus, we get
    $$t^{\frac{i}{2}}\|\bar{\nabla}^{i}\theta(t)\|_{C^{0}_{\bar{g}}}\longrightarrow 0,\quad \text{as}\quad t\longrightarrow 0,$$
    and
    $$t^{\frac{i+\alpha}{2}}[\bar{\nabla}^{i}\theta(t)]_{\alpha}\longrightarrow 0,\quad \text{as}\quad t\longrightarrow 0.$$
    Besides, we also get $\|\theta(t)\|_{C^{0}_{\bar{g}}}\longrightarrow\|\theta(0)\|_{C^{0}_{\bar{g}}}$ as $t\longrightarrow 0$. It follows that
    $$\lim_{t\rightarrow0}q_{n}(t)=\|\theta(0)\|_{C^{0}}.$$
    Now, we define $q_{n}(0)=\|\theta(0)\|_{C^{0}}$, one has
    $$\lim_{\tau\rightarrow0}K_{n}(\tau)=\|\theta(0)\|_{C^{0}}\leq \delta,$$
    which means $K_{n}(\tau)$ is a continuous function for $\tau$.

    Let $r$ be the radius of the neighborhood of $\bar{\varphi}$ preserved closedness in the moduli space of $G_{2}$-structures. Set $\kappa=\min\{1,\frac{r}{2}\}$. If $K_{n}(\tau)\leq \kappa$, then $K_{n}(\tau)\leq 1$ and 
    $$\|\theta(t)\|_{C^{0}}\leq K_{n}(\tau)\leq \frac{r}{2}<r,\quad t\in[0,\tau].$$
    Now, we choose $0<\delta\leq \delta_{\ast}$ and $\delta_{\ast}$ is defined by
    $$\delta_{\ast}:=\min\left\{\frac{\kappa}{4C_{0}},\frac{1}{8C_{0}C_{1}}\right\},$$
    where $C_{0}$ and $C_{1}$ are defined in Lemma \ref{lem4.1}. Now we prove \eqref{main} by contradiction. Assume there exists a time $\tau_{1}$ satisfying $0<\tau_{1}<T_{\max}$, $\tau_{1}\leq 1$ and
    \begin{align}
        K_{n}(\tau_{1})\geq 2C_{0}\delta.
    \end{align}
    We can define $\tau_{\ast}$ by
    $$\tau_{\ast}:=\inf\{t\in[0,\tau_{1}]\ |\ K_{n}(t)\geq 2C_{0}\delta\}.$$
    By the continuity of $K_{n}(t)$, we have
    $$K_{n}(\tau_{\ast})=2C_{0}\delta\leq 2C_{0}\delta_{\ast}\leq \frac{r}{2}<r.$$
    Lemma \ref{lem4.1} tells us that
    \begin{align}
        K_{n}(\tau_{\ast})&\leq C_{0}\delta+C_{1}K_{n}(\tau_{\ast})^{2}\notag\\
        &\leq C_{0}\delta+C_{1}(2C_{0}\delta)^{2}\notag\\
        &\leq C_{0}\delta+\frac{1}{2}C_{0}\delta\notag\\
        &=\frac{3}{2}C_{0}\delta\notag,
    \end{align}
    which contradicts $K_{n}(\tau_{\ast})=2C_{0}\delta$. Since $\tau_{1}$ is arbitrary, we prove \eqref{main}.

    It remains to prove $T_{\max}>1$. Suppose, for contradiction, that $T_{\max}\leq 1$, we consider $\frac{T_{\max}}{2}\leq t<T_{\max}$. By \eqref{main}, we have $K_{2}(t)\leq 2C_{0}\delta$, which means
    \begin{align}
        \sup_{\frac{T_{\max}}{2}\leq t<T_{\max}}\|\theta(t)\|_{C^{2,\alpha}}\leq CC_{0}\delta\left(1+T_{\max}^{-\frac{1}{2}}+T_{\max}^{-1}+T_{\max}^{-1-\frac{\alpha}{2}}\right)<+\infty.\notag
    \end{align}
    The last inequality uses $t\geq\frac{T_{\max}}{2}>0$. Besides, we also have
    $$\|\theta(t)\|_{C^{0}}\leq K_{n}(t)\leq 2C_{0}\delta<r.$$
    Let $0<\beta<\alpha$, and a sequence $t_{i}\rightarrow T_{\max}$ as $i\rightarrow\infty$. Since $M$ is compact, then  $C^{2,\alpha}$ is compactly embedded in $C^{2,\beta}$. Thus, we have
    $$\theta(t_{i})\longrightarrow\theta_{\ast},\quad \text{as} \quad i\longrightarrow \infty$$
    in the $C^{2,\beta}$ sense. Moreover, 
    $$d\theta_{\ast}=0, \quad \|\theta_{\ast}\|_{C^{0}}\leq 2C_{0}\delta<r.$$
    Since the equation \eqref{theta} is uniformly parabolic in the neighborhood $\|\theta\|_{C^0}<r$, local existence and continuous dependence provide a $C^{2,\beta}$-neighborhood $\mathcal{U}$ of $\theta_{\ast}$ and a common existence time $\tau>0$ for smooth closed initial values in $\mathcal{U}$. For sufficiently large $i$, we have
$$
\theta(t_i)\in\mathcal U,
\qquad
T_{\max}-t_i<\tau.
$$
Restarting the equation at $t_i$ and using uniqueness on the overlap extends the original solution beyond $T_{\max}$, a contradiction. Hence $T_{\max}>1$.
\end{proof}
Since Theorem \ref{thm1.1} has already been established above, it remains only to prove Theorem \ref{thm1.2}.

\begin{proof}[\bf Proof of Theorem \ref{thm1.2}]
Fix $\alpha\in(0,1)$ and apply Theorem \ref{thm1.1} with $n=9$. It follows that
\begin{equation}\label{eq:time-one}
\|\tilde{\varphi}(1)-\bar{\varphi}\|_{C^{9}}\leq C\|\varphi_{0}-\bar{\varphi}\|_{C^{0}}.
\end{equation}
Since the manifold $M$ is compact, we have
$$\|\tilde{\varphi}(1)-\bar{\varphi}\|_{W^{9,2}}
\le C_{\bar{g}}\|\tilde{\varphi}(1)-\bar{\varphi}\|_{C^{9}}.$$
By Theorem 1.3, Theorem 5.1 and Page 521-524 in \cite{LW19}, there is $\bar{\epsilon}>0$ such that a smooth structure in $[\bar{\varphi}]_+$ whose $C^9$ distance from $\bar{\varphi}$ is less than $\bar{\epsilon}$ generates a global solution of the Laplacian flow converging smoothly to $f^{\ast}\bar{\varphi}$ for some $f\in{\rm Diff}_0(M)$. 

Choose the $C^0$ neighborhood in the statement small enough that Theorem \ref{thm1.1} applies and the right side of \eqref{eq:time-one} is less than $\bar{\epsilon}$. Cohomology is preserved by \eqref{theta}. Hence there is a global solution $\hat{\varphi}(s)$ of the original Laplacian flow with
$$
\hat{\varphi}(0)=\tilde{\varphi}(1),\qquad
\hat{\varphi}(s)\longrightarrow f^{\ast}\bar{\varphi}\quad\text{smoothly as }s\to\infty.
$$
To recover the original flow before time one, solve the ordinary differential equation
\begin{equation}\label{eq:diffeomorphisms}
\frac{d}{dt}\chi_t=-V_{\bar{\varphi}}(\tilde{\varphi}(t))\circ\chi_t,
\qquad \chi_0={\rm id},
\end{equation}
on $0\leq t\leq 1+\varepsilon$, where $1+\varepsilon<T_{\max}$. The vector field is smooth on this compact time interval, so $\chi_t$ is a smooth family of diffeomorphisms in ${\rm Diff}_0(M)$. The pullback identity shows that
\[
\varphi(t)=\chi_t^{\ast}\tilde{\varphi}(t),\qquad 0\leq t\leq 1+\varepsilon,
\]
solves \eqref{eq:flow}. For $t\geq 1$, define
\[
\varphi(t)=\chi_1^{\ast}\hat{\varphi}(t-1).
\]
Diffeomorphism invariance makes this a solution of \eqref{eq:flow}, and its value at time one is the same as that of the first solution. Since $\tilde{\varphi}$ exists slightly beyond time one, uniqueness identifies the two solutions on a short overlap. They therefore give a smooth global solution starting from $\varphi_0$. Finally,
\[
\varphi(t)\longrightarrow\chi_1^{\ast}f^{\ast}\bar{\varphi}
=(f\circ\chi_1)^{\ast}\bar{\varphi}
\quad\text{smoothly},
\]
and $f\circ\chi_1\in{\rm Diff}_0(M)$. This proves the theorem.
\end{proof}

\textbf{Acknowledgments.}\ \ 
 The first author is supported by the China Postdoctoral Science Foundation (2026M793350). The second author is funded by the Shanghai Institute for Mathematics and Interdisciplinary Sciences (SIMIS) under grant number SIMIS-ID-2025-AD.


\begin{thebibliography}{99}

\bibitem{Bry06}
R. Bryant.
\textit{Some remarks on $G_{2}$-structures.}
Proceedings of G\"okova Geometry-Topology Conference 2005, 75–109, G\"okova Geometry/Topology Conference (GGT), G\"okova, 2006. 


\bibitem{BX11}
R. Bryant, F. Xu.
\textit{Laplacian flow for closed $G_{2}$-structures: short time behavior}.
arXiv:1101.2004 [math.DG].
 









\bibitem{FY18}
J. Fine, C. Yao.
\textit{Hypersymplectic $4$-manifolds, the $G_{2}$-Laplacian flow, and extension assuming bounded scalar curvature.}
Duke Math. J. {\bf167} (2018), no. 18, 3533–3589. 


\bibitem{Joy96}
D. D. Joyce.
\textit{Compact Riemannian $7$-manifolds with holonomy $G_{2}$. I, II.}
J. Differential Geom. {\bf 43} (1996), no. 2, 291–328, 329–375.



\bibitem{Hit00}
N. Hitchin.
\textit{The geometry of three-forms in six dimensions.}
J. Differential Geom. {\bf 55} (2000), no. 3, 547–576. 
 
 

\bibitem{Spiros introduce to G2}
S. Karigiannis.
\textit{Introduction to $G_{2}$ geometry.}
Lectures and surveys on $G_{2}$-manifolds and related topics, 3–50, Fields Inst. Commun., {\bf 84}, Springer, New York, 2020. 





\bibitem{LL21}
B. Lambert, J. Lotay. 
\textit{Spacelike mean curvature flow.}
J. Geom. Anal. {\bf 31} (2021), no. 2, 1291–1359. 




\bibitem{LL26}
C. Li, Y. Li.
\emph{Curvature pinching estimate under the Laplacian $G_2$ flow.}
J. Funct. Anal. \textbf{290} (2026), Article 111199.




\bibitem{LLZ26}
C. Li, Y. Li, K. Zheng.
\emph{Long time existence of the Laplacian flow with Weyl tensor.}
arXiv:2610.23412 [math.DG]


\bibitem{local curvature of G2}
Y. Li.
\textit{Local curvature estimates for the Laplacian flow.} 
Calc. Var. Partial Differential Equations {\bf 60} (2021), no. 1, Paper No. 28, 37 pp. 






\bibitem{LW17}
J. Lotay, Y. Wei.
\emph{Laplacian flow for closed $G_2$ structures: Shi-type estimates,
uniqueness and compactness.}
Geom. Funct. Anal. \textbf{27} (2017), no. 1, 165-233.


\bibitem{LW19}
J. Lotay, Y. Wei. 
\textit{Stability of torsion-free $G_{2}$ structures along the Laplacian flow.}
 J. Differential Geom. {\bf 111} (2019), no. 3, 495–526.


\bibitem{PS24}
S. Picard, C. Suan.
 \textit{Flows of $G_{2}$-Structures associated with Calabi-Yau Manifolds.}
 Math. Res. Lett. {\bf 31} (2024), no. 6, 1837–1877.

\bibitem{Ros97}
S. Rosenberg.
\textit{The Laplacian on a Riemannian Manifold: An Introduction to Analysis on Manifolds.}
Cambridge University Press, 1997.

\end{thebibliography}
\end{document}